\documentclass[11pt]{amsart}

\usepackage{amsmath,amsthm, amssymb, amsfonts}
\usepackage[mathscr]{eucal}
\usepackage{epsfig}
\usepackage{graphicx}
\usepackage{psfrag}
\usepackage[usenames]{color}
\usepackage{subfigure}

\input epsf

\usepackage[colorlinks=true, urlcolor=NavyBlue, linkcolor=NavyBlue, citecolor=NavyBlue, pdftitle={On the Transverse Invariant for Bindings of Open Books}, pdfauthor={David Shea Vela--Vick}, pdfsubject={Transverse Invariants}, pdfkeywords={Legendrian knots, transverse knots, Heegaard Floer homology, 57M27, 57R58}]{hyperref}

\newtheorem{theorem}{Theorem}
\newtheorem{corollary}[theorem]{Corollary}

\newtheorem{othertheorem}{Theorem}[section]

\newtheorem{othercorollary}[othertheorem]{Corollary}
\newtheorem{lemma}[othertheorem]{Lemma}

\theoremstyle{definition}
\newtheorem{definition}{Definition}[section]

\theoremstyle{definition}
\newtheorem{remark}[othertheorem]{Remark}

\newcommand{\thmref}[1]{Theorem~\ref{#1}}

\newcommand{\mfr}[1]{\mathfrak{#1}}
\newcommand{\mr}[1]{\mathrm{#1}}

\newcommand{\HFKH}{\widehat{\mathrm{HFK}}}
\newcommand{\HFKM}{\mathrm{HFK}^-}

\newcommand{\HFH}{\widehat{\mathrm{HF}}}
\newcommand{\STH}{\widehat{\mathscr{T}}}
\newcommand{\ST}{\mathscr{T}}
\newcommand{\SLH}{\widehat{\mathscr{L}}}
\newcommand{\SL}{\mathscr{L}}

\begin{document}

\thispagestyle{empty}

\title[On the Transverse Invariant]{On the Transverse Invariant for Bindings of Open Books}
\author{David Shea Vela--Vick}
\address{Department of Mathematics \\ Columbia University}
\email{shea@math.columbia.edu}
\urladdr{\href{http://www.math.columbia.edu/~shea}{http://www.math.columbia.edu/\~{}shea}}

\date{\today}
\maketitle

\begin{abstract}

Let $T \subset (Y,\xi)$ be a transverse knot which is the binding of some open book, $(T,\pi)$, for the ambient contact manifold $(Y,\xi)$.  In this paper, we show that the transverse invariant $\STH(T) \in \HFKH(-Y,K)$, defined in \cite{LOSS}, is nonvanishing for such transverse knots.  This is true regardless of whether or not $\xi$ is tight.  We also prove a vanishing theorem for the invariants $\SL$ and $\ST$.  As a corollary, we show that if $(T,\pi)$ is an open book with connected binding, then the complement of $T$ has no Giroux torsion.

\end{abstract}

\section{Introduction} 
	\label{sec:intro}

In a recent paper by Lisca, Ozsv\'ath, Stipsicz, and Szab\'o \cite{LOSS}, the authors define invariants of null-homologous Legendrian and transverse knots.  These invariants live in the knot Floer homology groups of the ambient space with reversed orientation, and generalize the previously defined invariants of closed contact manifolds, $c(Y,\xi)$.  They have been useful in constructing new examples of knot types which are not transversally simple (see \cite{LOSS,OSt}), and play an important role in the classification of Legendrian and transverse twist knots (see \cite{ENV}).

In this paper, we investigate properties of these invariants for a certain important class of transverse knots.

\begin{theorem}\label{thm:binding}
Let $T \subset (Y,\xi)$ be a transverse knot which can be realized as the binding of an open book $(T,\pi)$ compatible with the contact structure $\xi$.  Then, the transverse invariant $\STH(T)$ is nonvanishing.
\end{theorem}

\begin{remark}
In \cite{LOSS}, it is shown that if $c(Y,\xi) \neq 0$, then $\ST(T) \neq 0$ for any transverse knot $T \subset (Y,\xi)$.  In \thmref{thm:binding}, no restrictions are made on the ambient contact structure $\xi$.  In particular, the theorem is true even when $\xi$ is overtwisted.  Moreover, the nonvanishing of the invariant $\STH$ implies the nonvanishing of the invariant $\ST$.
\end{remark}

Let $L$ be a null-homologous Legendrian knot in $(Y,\xi)$.  It is shown in \cite{LOSS} that the invariant $\SL(L)$ inside $\HFKM(-Y,L)$ remains unchanged under negative stabilization, and therefore yields an invariant of transverse knots.  If $T$ is a null-homologous transverse knot in $(Y,\xi)$ and $L$ is a Legendrian approximation of $T$, then $\ST(T) := \SL(L)$.  We will generally state results only in the Legendrian case, even though the same results are also true in the transverse case.

\begin{remark}
There is a natural map $\HFKM(-Y,L) \to \HFKH(-Y,L)$, induced by setting $U = 0$.  Under this map, the $\SL(L)$ is sent to $\SLH(L)$.  Therefore, if $\SLH(L)$ is nonzero, then $\SL(L)$ must also be nonzero.  Similarly, $\SL(L)$ vanishing implies that $\SLH(L)$ must also vanish.
\end{remark}

In addition to understanding when these invariants are nonzero, we are also interested in circumstances under which they vanish.  In \cite{LOSS}, it was shown that if the complement of a Legendrian knot contains an overtwisted disk, then the Legendrian invariant for that knot vanishes.  Here, we generalize this result by proving:

\begin{theorem}\label{thm:vanishing}
Let $L$ be a Legendrian knot in a contact manifold $(Y,\xi)$.  If the complement $Y-L$ contains a compact submanifold $N$ with convex boundary such that $c(N,\xi|_N) = 0$ in $SFH(-N,\Gamma)$, then the Legendrian invariant $\SL(L)$ vanishes.
\end{theorem}

Since $I$-invariant neighborhoods of convex overtwisted disks have vanishing contact invariant (Example 1 of \cite{HKM2}), \thmref{thm:vanishing} generalizes the vanishing theorem from \cite{LOSS}.

In \cite{GHV}, Ghiggini, Honda, and Van Horn--Morris show that a closed contact manifold with positive Giroux torsion has vanishing contact invariant.  They show this by proving that the contact element for a $2\pi$-torsion layer vanishes in sutured Floer homology.  Thus, as an immediate corollary to \thmref{thm:vanishing}, we have:

\begin{corollary}\label{cor:torsion}
Let $L$ be a Legendrian knot in a contact manifold $(Y,\xi)$.  If the complement $Y-L$ has positive Giroux torsion, then the Legendrian invariant $\SL(L)$ vanishes.
\end{corollary}

\begin{remark}
A similar result has been established for the weaker invariant $\SLH$ by Stipsicz and V\'ertesi \cite{StV} using slightly different arguments.
\end{remark}

\begin{remark}
\thmref{thm:binding} and Corollary \ref{cor:torsion} are both true in the transverse case as well.
\end{remark}

Combining the transverse version of Corollary \ref{cor:torsion} with \thmref{thm:binding}, we conclude the following interesting fact about complements of connected bindings:

\begin{theorem}\label{thm:torsion2}
Let $(T,\pi)$ be an open book with a single binding component.  Then the complement of $T$ has no Giroux torsion.
\end{theorem}

As Giroux torsion is presently the only known mechanism for a 3-manifold to admit more than a finite number of tight contact structures, it is important to understand the relationship between tight contact structures with positive Giroux torsion and the open books which support them.  Of course, Theorem~\ref{thm:torsion2} only applies to connected bindings of open books, leading one to conjecture that it should be true for arbitrary open book decompositions.  We prove this with Etnyre in \cite{EV} using different methods.

\begin{othertheorem}[Etnyre-Vela--Vick, \cite{EV}]\label{thm:torsion3}
Let $(B,\pi)$ be an open book for a contact manifold $(Y,\xi)$.  Then the complement of $B$ has no Giroux torsion.
\end{othertheorem}

This paper is organized as follows: in Section \ref{sec:background}, we briefly review some of the basic concepts in contact geometry and knot Floer homology.  Section \ref{sec:proof_of_main_thm} is devoted to proving \thmref{thm:binding}.  In Section \ref{sec:vanishing}, we conclude with a proof of \thmref{thm:vanishing}.

\section*{Acknowledgements} 
\label{sec:acknowledgements}
I owe a tremendous debt of gratitude to my advisor, John Etnyre.  This problem arose in discussions with John.  His support and guidance over the years have been warmly received and much appreciated.  I also thank Clayton Shonkwiler for providing valuable comments on drafts of this paper.


\section{Background} 
\label{sec:background}

\subsection{Contact Geometry Preliminaries} 
\label{sub:contact_geometry_preliminaries}

Recall that a contact structure on an oriented 3-manifold is a plane field $\xi$ satisfying a certain nonintegrability condition.  We assume that our plane fields are cooriented, and that $\xi$ is given as the kernel of some global 1-form: $\xi = \mr{ker}(\alpha)$ with $\alpha(N_p) > 0$ for each oriented normal vector $N_p$ to $\xi_p$.  Such an $\alpha$ is called a {\it contact form for $\xi$}.  In this case, the nonintegrability condition is equivalent to the statement $\alpha \wedge d\alpha >0$.

A primary tool used in the study of contact manifolds has been Giroux's correspondence between contact structures on 3-manifolds and open book decompositions up to an equivalence called {\it positive stabilization}.  An {\it open book decomposition} of a contact 3-manifold $(Y,\xi)$ is a pair $(B,\pi)$, where $B$ is an oriented, fibered, transverse link and $\pi : (Y - B) \to S^1$ is a fibration of the complement of $B$ by oriented surfaces whose oriented boundary is $B$.  

An open book is said to be {\it compatible} with a contact structure $\xi$ if, after an isotopy of the contact structure, there exists a contact form $\alpha$ for $\xi$ such that:
\begin{enumerate}
	\item $\alpha(v)>0$ for each (nonzero) oriented tangent vector $v$ to $B$, and 
	\item $d\alpha$ restricts to a positive area form on each page of the open book.  
\end{enumerate}

Given an open book decompositon of a 3-manifold $Y$, Thurston and Winkelnkemper \cite{TW} show how one can produce a compatible contact structure on $Y$.  Giroux proves in \cite{Gi} that two contact structures which are compatible with the same open book are, in fact, isotopic as contact structures.  Giroux also proves that two contact structures are isotopic if and only if they are compatible with open books which are related by a sequence of positive stabilizations.

\begin{definition}\label{def:Legendrian}
	A knot $L$ smoothly embedded in a contact manifold $(Y,\xi)$ is called {\it Legendrian} if $T_pL \subset \xi_p$ for all $p$ in $L$.
\end{definition}

\begin{definition}\label{def:transverse}
	An oriented knot $T$ smoothly embedded in a contact manifold $(Y,\xi)$ is called {\it transverse} if it always intersects the contact planes transversally with each intersection positive.
\end{definition}

We say that two Legendrian knots are {\it Legendrian isotopic} if they are isotopic through Legendrian knots; similarly, two transverse knots are {\it transversally isotopic} if they are isotopic through transverse knots.  Given a Legendrian knot $L$, one can produce a canonical transverse knot nearby to $L$, called the {\it transverse pushoff} of $L$.  On the other hand, given a transverse knot $T$, there are many nearby Legendrian knots, called {\it Legendrian approximations} of $T$.  Although there are infinitely many distinct Legendrian approximations of a given transverse knot, they are all related to one another by sequences of negative stabilizations.  These two constructions are inverses to one another, up to the ambiguity involved in choosing a Legendrian approximation of a given transverse knot (see \cite{EFM, EH}).

If $I$ is an invariant of Legendrian knots which remains unchanged under negative stabilization, then $I$ is also an invariant of transverse knots: if $T$ is a transverse knot and $L$ is a Legendrian approximation of $T$, define $I(T)$ to be equal to the invariant $I(L)$ of the Legendrian knot $L$.  This is how the authors define the transverse invariants $\ST(T)$ and $\STH(T)$ in \cite{LOSS}.  

For more on open book decompositions and on Legendrian and transverse knots, we refer the reader to \cite{Et1,Co} and, respectively, \cite{Et2}.  


\subsection{Heegaard Floer Preliminaries} 
\label{sub:heegaard_floer_preliminaries}

This paper is primarily concerned with two versions of Heegaard Floer homology, which are invariants of (null-homologous) knots inside closed 3-manifolds.  These homologies, called knot Floer homology, are denoted $\HFKM(Y,K)$ and $\HFKH(Y,K)$.  In knot Floer homology, the basic input is a doubly-pointed Heegaard diagram; that is, a Heegaard diagram $(\Sigma,\alpha,\beta)$, together with two basepoints $w, z \in \Sigma - (\alpha \cup \beta)$, in the complement of the $\alpha$- and $\beta$-curves.  These diagrams are required to satisfy certain admissibility conditions which depend on the version of the theory which one is working.

Given a doubly pointed Heegaard diagram, one can produce a knot in the resulting 3-manifold.  To do this, connect $z$ to $w$ by an arc in the complement of the $\alpha$-curves, and $w$ to $z$ by an arc in the complement of the $\beta$-curves.  After depressing the interiors of these arcs into the $\beta$- and $\alpha$-handlebodies, respectively, the result is an oriented knot inside the closed 3-manifold specified by the Heegaard diagram $(\Sigma,\alpha,\beta)$.  Using a bit of elementary Morse theory, one can show that any knot in any closed 3-manifold can be represented by a doubly-pointed Heegaard diagram.

If the genus of $\Sigma$ is $g$, then the chain groups for $\HFKM(Y,K)$ are generated as a $\mathbb{Z}/2[U]$-module by the intersection points between the two $g$-dimensional subtori $\mathbb{T}_\alpha = \alpha_1 \times \dots \times \alpha_g$ and $\mathbb{T}_\beta = \beta_1 \times \dots \times \beta_g$ inside $\mr{Sym}^g(\Sigma)$.  Given a complex structure on $\Sigma$, $\mr{Sym}^g(\Sigma)$ inherits a natural complex structure from the projection $\times_g \Sigma \to \mr{Sym}^g(\Sigma)$.  The boundary map counts certain rigid holomorphic disks in $\mr{Sym}^g(\Sigma)$, with boundary lying on $\mathbb{T}_\alpha \cup \mathbb{T}_\beta$, connecting these intersection points: 
\[
	\partial^- {\bf x} = 
	\sum_{{\bf y} \in \mathbb{T}_\alpha \cap \mathbb{T}_\beta} 
	\sum_{\substack{\phi \in \pi_2({\bf x}, {\bf y}), \,
	\mu(\phi)=1, \\
	n_z(\phi) = 0}} 
	\# \widehat{\mfr{M}}(\phi) \cdot U^{n_w(\phi)} \cdot {\bf y}.
\]
Here $n_v(\phi)$ is equal to the algebraic number of times the disk $\phi$ intersects the subspace $\{v\} \times \mr{Sym}^{g-1}(\Sigma)$; $\pi_2({\bf x}, {\bf y})$ is the set of homotopy classes of disks connecting ${\bf x}$ to ${\bf y}$ with boundaries lying on $\mathbb{T}_\alpha$ and $\mathbb{T}_\beta$. 

The chain groups for $\HFKH(Y,K)$ are generated as a $\mathbb{Z}/2$-vector space by the intersection points between $\mathbb{T}_\alpha$ and $\mathbb{T}_\beta$ in $\mr{Sym}^g(\Sigma)$.  In this case, the boundary map counts holomorphic disks in $\mr{Sym}^g(\Sigma)$, with boundaries lying on $\mathbb{T}_\alpha \cap \mathbb{T}_\beta$, missing both $z$ and $w$:
\[
	\widehat{\partial} {\bf x} = 
	\sum_{{\bf y} \in \mathbb{T}_\alpha \cap \mathbb{T}_\beta} 
	\sum_{\substack{\phi \in \pi_2({\bf x}, {\bf y}), \,
	\mu(\phi) = 1, \\
	n_z(\phi) = 0, \,
	n_w(\phi) = 0}} 
	\# \widehat{\mfr{M}}(\phi) \cdot {\bf y}.
\]

For more information on Heegaard Floer homology and knot Floer homology, we refer the reader to \cite{OS1,Li} and \cite{OS3, Ra}, respectively.


\subsection{Invariants of Legendrian and Transverse Knots} 
\label{sub:invariants_of_legendrian_and_transverse_knots}

Let $L$ be a Legendrian knot with knot type $K$, and let $T$ be a transverse knot in the same knot type.  In \cite{LOSS}, the authors define invariants $\SL(L)$ and $\ST(T)$ in $\HFKM(-Y,K)$ and $\SLH(L)$ and $\STH(T)$ in $\HFKH(-Y,K)$.  These invariants are constructed in a similar fashion to the contact invariants in \cite{HKM1,HKM2}.  Below we describe how to construct the invariant for a Legendrian knot.

Let $L \subset (Y,\xi)$ be a null-homologous Legendrian knot.  Consider an open book decomposition of $(Y,\xi)$ containing $L$ on a page $S$.  Choose a basis $\{a_0,\dots,a_n\}$ for $S$ (i.e a collection of disjoint, properly embedded arcs $\{a_0, \dots, a_n \}$ such that $S - \bigcup a_i$ is homeomorphic to a disk) with the property that $L$ intersects only the basis element $a_0$, and does so transversally in a single point.  Let $\{b_0,\dots,b_n\}$ be a collection of properly embedded arcs obtained from the $a_i$ by applying a small isotopy so that the endpoints of the arcs are isotoped according to the induced orientation on $\partial S$ and so that each $b_i$ intersects $a_i$ transversally in the single point $x_i$.  If $\phi: S \to S$ is the monodromy map representing the chosen open book decomposition, then our Heegaard diagram is given by
\[
	(\Sigma,\alpha,\beta) = (S_{1/2} \cup -S_1, (a_i \cup a_i), (b_i \cup \phi(b_i))).
\]

The first basepoint, $z$, is placed on the page $S_{1/2}$ in the complement of the thin strips of isotopy between the $a_i$ and $b_i$.  The second basepoint, $w$, is placed on the page $S_{1/2}$ inside the thin strip of isotopy between $a_0$ and $b_0$.  The two possible placements of $w$ correspond to the two possible orientations of $L$.

The Lengendrian invariant $\SL(L)$ is defined, up to isomorphism, to be the element $[{\bf x}] = [(x_0, \dots, x_n)]$ in $\HFKM(\Sigma,\beta,\alpha,w,z)$.  A picture of this construction in the case at hand is given in Figure~\ref{fig:basis}.  If $T$ is a transverse knot, the transverse invariant $\ST(T)$ is defined to be the Legendrian invariant of a Legendrian approximation of $T$.

One interesting property of these invariants is that they do not necessarily vanish for knots in an overtwisted contact manifolds; this is why we do not need to assume tightness in \thmref{thm:binding}.  Another property, which will be useful in Section \ref{sec:proof_of_main_thm}, is that these invariants are natural with respect to contact $(+1)$-surgeries.

\begin{othertheorem}[Ozsv\'ath-Stipsicz, \cite{OSt}]\label{thm:surgery}
Let $L \subset (Y,\xi)$ be a Legendrian knot.  If $(Y_B, \xi_B, L_B)$ is obtained from $(Y,\xi,L)$ by contact ($+1$)-surgery along a Legendrian knot $B$ in $(Y,\xi,L)$, then under the natural map
\[
	F_B : \HFKM(-Y,L) \to \HFKM(-Y_B,L_B),
\]
$\SL(L)$ is mapped to $\SL(L_B)$.
\end{othertheorem}

An immediate corollary to this fact is the following:

\begin{othercorollary}\label{cor:legsurgery}
Let $L \subset (Y,\xi)$ be a Legendrian knot, and suppose that $(Y_B, \xi_B, L_B)$ is obtained from $(Y,\xi,L)$ by Legendrian surgery along a Legendrian knot $B$ in $(Y,\xi,L)$.  If $\SL(L) \neq 0$ in $\HFKM(-Y,L)$, then $\SL(L_B) \neq 0$ in $\HFKM(-Y_B,L_B)$.
\end{othercorollary}

\begin{remark}
\thmref{thm:surgery} and Corollary \ref{cor:legsurgery} are also true for the invariant $\SLH(L)$ and for the invariants $\ST(T)$ and $\STH(T)$ in the case of a transverse knot.
\end{remark}

In addition, the invariant $\SL$ directly generalizes the original contact invariant $c(Y,\xi) \in \HFH(-Y)$ (see \cite{OS4}).  Under the natural map $\HFKM(-Y,L) \to \HFH(-Y)$ induced by setting $U = 1$, $\SL(L)$ maps to $c(Y,\xi)$, the contact invariant of the ambient contact manifold.

We encourage the interested reader to look at \cite{LOSS,OSt} to learn about other properties of these invariants.


\section{Proof of \thmref{thm:binding}} 
\label{sec:proof_of_main_thm}

Let $T \subset (Y,\xi)$ be a transverse knot.  Recall that \thmref{thm:binding} states that if $T$ is the binding for some open book $(T,\pi)$ for $(Y,\xi)$, then the transverse invariant $\STH(T) \in \HFKH(-Y,T)$ is nonvanishing.  

In this section, we prove \thmref{thm:binding} in three steps.  In Section \ref{sub:obtaining_the_pointed_diagram} we construct an open book on which a Legendrian approximation $L$ of the transverse knot $T$ sits.  Then we show in Section \ref{sub:admissibility} that the Heegaard diagram obtained in Section \ref{sub:obtaining_the_pointed_diagram} is weakly admissible.  Finally, in Section \ref{sub:computing_T(T)}, we prove the theorem in the special case where the monodromy map $\phi_n$ consists of a product of $n$ negative Dehn twists along a boundary-parallel curve.

An arbitrary monodromy map differs from some such $\phi_n$ by a sequence of positive Dehn twists, or Legendrian surgeries, along curves contained in pages of the open book.  By Corollary \ref{cor:legsurgery}, since the transverse invariant is nonvanishing for the monodromy maps $\phi_n$, it must also be nonvanishing for an arbitrary monodromy map.

\subsection{Obtaining the pointed diagram} 
\label{sub:obtaining_the_pointed_diagram}

By hypothesis, $T$ is the binding of an open book $(T, \pi)$ for $(Y,\xi)$.  To compute the transverse invariant $\STH(T)$, we need to find a Legendrian approximation $L$ of $T$, realized as a curve on a page of an open book for $(Y,\xi)$.

\begin{figure}[htbp]
	\centering
	\subfigure[ ]{\label{fig:OpenBook1} 
	\psfrag{a}{\hspace{-5px}$T$}
	\psfrag{b}{$\gamma$}
	\psfrag{c}{$S$}
	\includegraphics[scale=0.55]{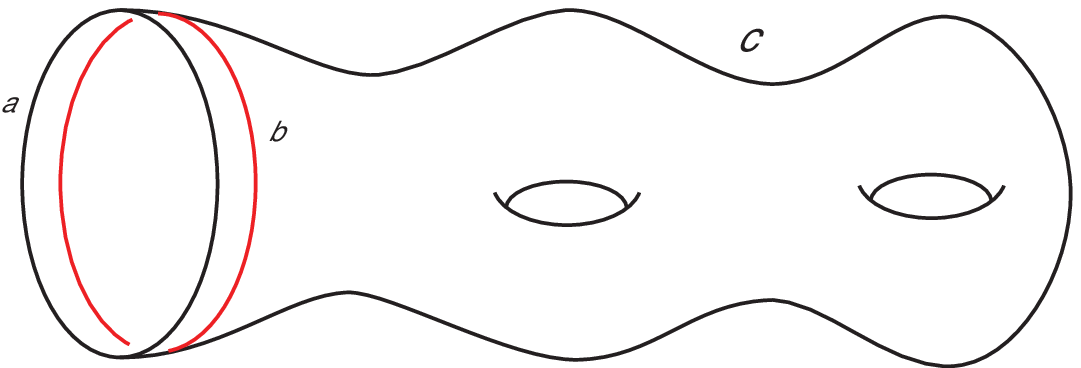}}
	\subfigure[ ]{\label{fig:OpenBook2}
	\psfrag{a}{$\gamma$}
	\psfrag{b}{$c$}
	\psfrag{c}{$\tau$}
	\psfrag{d}{$T$}
	\psfrag{e}{$S^\prime$}
	\includegraphics[scale=0.55]{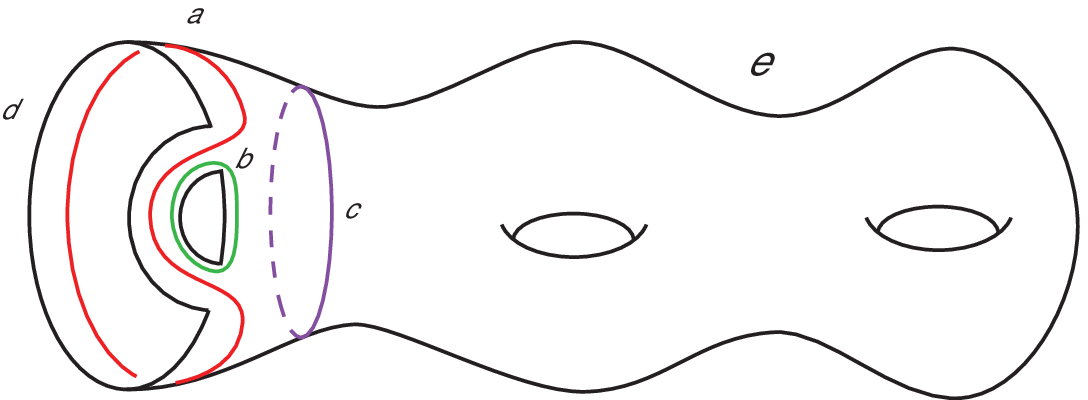}}
	\caption{ }
	\label{fig:OpenBook}
\end{figure}

In Figure~\ref{fig:OpenBook1}, we see a page of the open book $(T,\pi)$.  Here, $T$ appears as the binding $\partial S = T$.  Assuming the curve $\gamma$ could be realized as a Legendrian curve, it would be the natural choice for the Legendrian approximation $L$.  Unfortunately, since $\gamma$ is zero in the homology of the page, $\gamma$ cannot be made Legendrian on the page.  

To fix this problem, stabilize the diagram.  The result of such a stabilization is illustrated in Figure~\ref{fig:OpenBook2}.  To see that this solves the problem, we prove the following lemma:

\begin{lemma}\label{lem:stabilization}
	The stabilization depicted in Figure~\ref{fig:OpenBook2} can be performed while fixing $T$ as the ``outer'' boundary component.  
\end{lemma}

Assume the truth of Lemma \ref{lem:stabilization} for the moment.  Then the curve $\gamma$ depicted in Figure~\ref{fig:OpenBook2} can now be Legendrian realized, as it now represents a nonzero element in the homology of the page.  By construction, if we orient this Legendrian coherently with $T$, then $T$ is the transverse pushoff of $\gamma$.

\begin{proof}
Consider $S^3$ with its standard tight contact structure.  Let $(H_+,\pi_+)$ be the open book for $(S^3,\xi_{\mr{std}})$ whose binding consists of two perpendicular Hopf circles and whose pages consist of negative Hopf bands connecting these two curves.  In this case, each binding component is a transverse unknot with self-linking number equal to $-1$.

Let $T$ be a transverse knot contained in a contact manifold $(Y,\xi)$ and let $U$ be a transverse unknot in $(S^3,\xi|_{\mr{std}})$ with self-linking number equal to $-1$.  Observe that the complement of a standard neighborhood of a point contained in $U$ is itself a standard neighborhood of a point contained in a transverse knot.  Therefore, if we perform a transverse connected sum of $T$ with the transverse unknot of self-linking number equal to $-1$ in $(S^3, \xi_{\mr{std}})$, we do not change the transverse knot type of $T$.  

\begin{figure}[htbp]
	\centering
	\subfigure[ ]{\label{fig:murasugi1} 
		\hspace{-6pt}\includegraphics[scale=1.15]{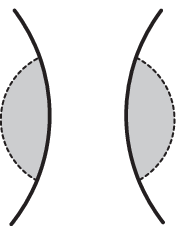}}
	\hspace{2cm}
	\subfigure[ ]{\label{fig:murasugi2}
		\hspace{-6pt}\includegraphics[scale=1.15]{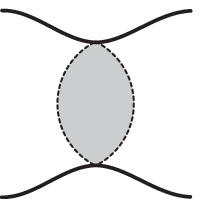}}
	\caption{ }
	\label{fig:Murasugi}
\end{figure}

Let $(B,\pi)$ be an open book with connected binding for a contact manifold $(Y,\xi)$.  Consider the contact manifold obtained from $(Y,\xi)$ by Murasugi summing the open book $(B,\pi)$ with the open book $(H_+,\pi_+)$ along bigon regions bounded by boundary-parallel arcs contained in pages of the respective open books.  The summing process is depicted in Figure~\ref{fig:Murasugi}.  Figure~\ref{fig:murasugi1} shows the open books before the Murasugi sum, while Figure~\ref{fig:murasugi2} shows the resulting open book after the sum. 

\begin{figure}[htbp]
	\centering
	\subfigure[ ]{\label{fig:binding1} 
		\hspace{-7pt}\includegraphics[scale=0.6]{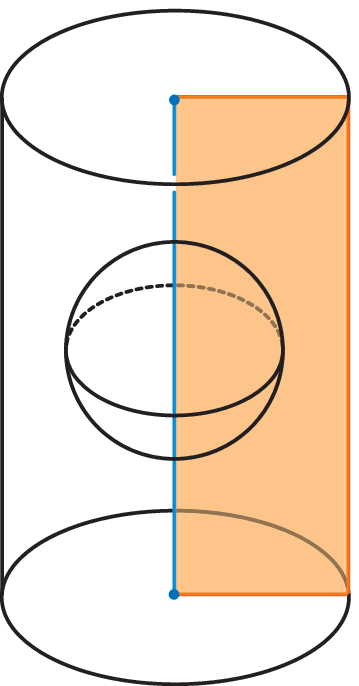}}
	\hspace{2cm}
	\subfigure[ ]{\label{fig:binding2}
		\hspace{-7pt}\includegraphics[scale=0.6]{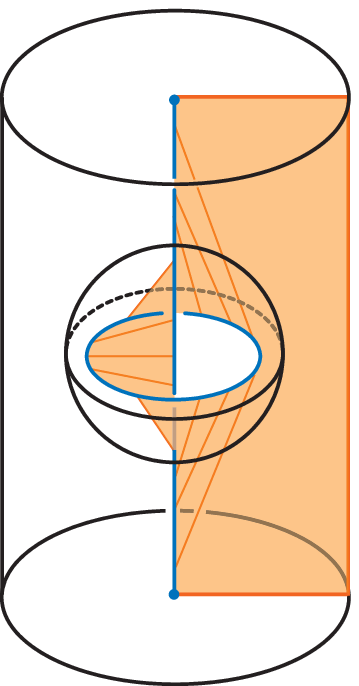}}
	\caption{ }
	\label{fig:binding}
\end{figure}

The Murasugi sum operation has the effect of performing a contact connected sum of $(Y,\xi)$ with $(S^3,\xi_{\mr{std}})$ and a transverse connected sum of the binding component $B$ with one of the binding components of $H_+$ (see \cite{To}).  Before and after pictures of this operation are shown in Figures \ref{fig:binding1} and \ref{fig:binding2}, respectively.  Since both of the binding components of the open book $(H_+,\pi_+)$ are transverse unknots with self-linking number equal to $-1$, this connected sum has no effect on the transverse knot type of the ``outer'' boundary component of the open book in Figure~\ref{fig:OpenBook}.
\end{proof}

Since the curve $\gamma$ can now be Legendrian realized and approximates $T$ as desired, we will denote $\gamma$ by $L$ from this point forward.  The new monodromy map $\phi^\prime : S^\prime \to S^\prime$ is equal to the old monodromy map, $\phi$, composed with one positive Dehn twist along the curve $c$ shown in Figure~\ref{fig:OpenBook2}.  For notational ease, we continue denoting the monodromy map by $\phi$, and the page by $S$.

\begin{figure}[htbp]
	\centering
	\psfrag{a}{$a_o$}
	\psfrag{b}{$w_1$}
	\psfrag{c}{$w_2$}
	\psfrag{d}{$z$}
	\psfrag{e}{$a_{2g}$}
	\psfrag{f}{$x_{2g}$}
	\psfrag{g}{$a_1$}
	\psfrag{h}{$x_1$}
	\psfrag{i}{$\dots$}
	\psfrag{j}{$L$}
	\psfrag{k}{$T$}
	\psfrag{l}{$x_0$}
        \includegraphics[scale=1.25]{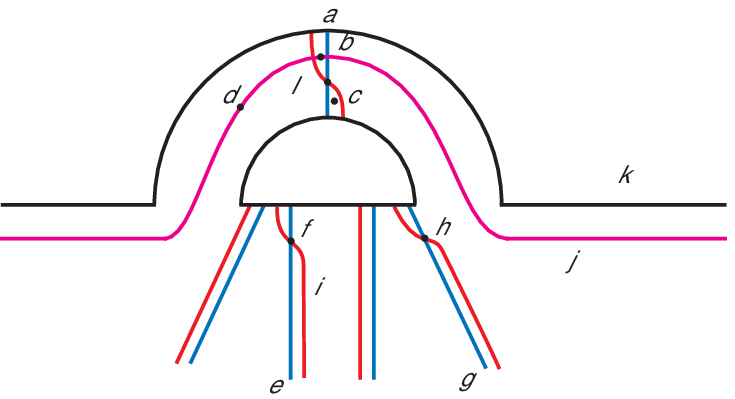}
		\caption{}
		\label{fig:basis}
\end{figure}

We choose a basis for our surface whose local picture near the stabilization is depicted in Figure~\ref{fig:basis}.  There are two possible choices for the placement of the second basepoint $w$: $w_1$ and $w_2$.  In order for $L$ to be oriented coherently with $T$, we must choose $w = w_1$.


\subsection{Admissibility} 
\label{sub:admissibility}

Our goal is to construct a weakly admissible, doubly-pointed Heegaard diagram from the open book described in Section \ref{sub:obtaining_the_pointed_diagram}.  

Before we continue, let us discuss some notation.  We are concerned with open book decompositions whose pages are twice-punctured surfaces.  We picture a genus $g$ surface as a $4g$-sided polygon with certain boundary edges identified.  We choose the standard identification scheme, where the first and third edges are identified, as are the second and fourth edges, the fifth and seventh edges, and so on.  For convenience, we always assume that the first edge appears in the 12 o'clock position, at the top of each diagram.

\begin{figure}[htbp]
	\centering
	\psfrag{a}{$x_1$}
	\psfrag{b}{$x_2$}
	\psfrag{c}{$z$}
	\psfrag{d}{$x_0$}
	\psfrag{e}{$w$}
	\psfrag{f}{{\Large $S_{1/2}$}}
	\psfrag{g}{$a_{1,I}$}
	\psfrag{h}{$a_{2,I}$}
	\psfrag{i}{$a_{1,F}$}
	\psfrag{j}{$a_{2,F}$}
	\psfrag{k}{$a_0$}
	\psfrag{l}{$b_0$}
	    \includegraphics[scale=1.2]{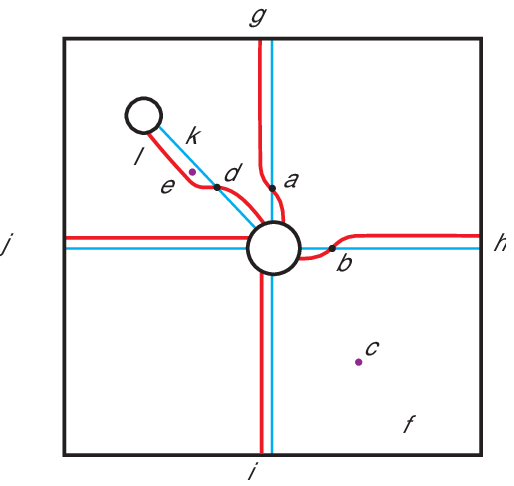}
		\caption{}
		\label{fig:S_one_half}
\end{figure}

Our punctures are always situated so that one of the punctures is in the center of the polygon, with the other close by.  We choose our basis elements, $a_1, \dots, a_{2g}$, to be straight arcs emanating from the center of the polygon and passing out the corresponding edge.  If we were to forget about the identifications being made at the boundary, the basis element $a_i$ would break into two straight arcs emanating from the center of the diagram.  For ease of exposition in what follows, we label the first segment that we see as we move clockwise around the diagram $a_{i,I}$, and the second $a_{i,F}$, where the subscript $I$ stands for ``initial'', while the subscript $F$ stands for ``final''.  

Up to isotopy, we may assume that the second boundary component of our surface lies (pictorially) in the region between the curves $a_{2g,F}$ and $a_{1,I}$, as shown in Figure~\ref{fig:S_one_half}.  The last basis element $a_0$ is a straight line segment connecting the two boundary components of the surface.

We have adopted the practice of Honda, Kazez and Mati\'c of placing surrogate basepoints throughout the diagram whenever it is convenient.  This signals that the local multiplicity of any domain contributing to the differential is zero in that region.

\begin{figure}[htbp]
	\centering
	\psfrag{a}{$a_{1,I}$}
	\psfrag{b}{$a_{2,I}$}
	\psfrag{c}{$a_{1,F}$}
	\psfrag{d}{$a_{2,F}$}
	\psfrag{e}{$(a_{1,I})$}
	\psfrag{f}{$(a_{2,I})$}
	\psfrag{g}{$(a_{1,F})$}
	\psfrag{h}{$(a_{2,F})$}
	\psfrag{i}{{\Large $-S_1$}}
	\psfrag{j}{$U$}
	\psfrag{k}{$R$}
	\psfrag{l}{$\gamma$}
	\psfrag{m}{$\delta$}
        \includegraphics[scale=1.65]{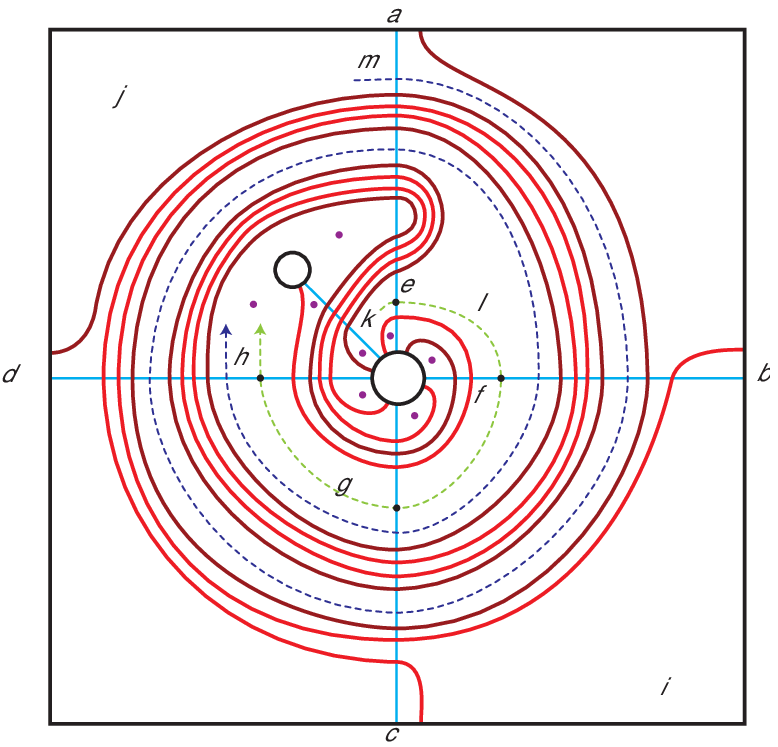}
		\caption{}
		\label{fig:BigPicture}
\end{figure}

We have restricted our figures to the case where our page is a twice-punctured torus, and our monodromy map $\phi$ consists of two negative Dehn twists along the curve $\tau$ in Figure~\ref{fig:OpenBook2}.  The resulting doubly-pointed Heegaard diagram is shown in Figures \ref{fig:S_one_half} and \ref{fig:BigPicture}.  

Figure~\ref{fig:S_one_half} shows the $S_{1/2}$ page of our open book, while Figure~\ref{fig:BigPicture} shows the $-S_1$ page (note the reversed orientation).  The invariant appears in Figure~\ref{fig:S_one_half} as the intersection point ${\bf x} = (x_0,x_1,x_2)$.

Consider the small region southeast of $x_0$ in Figure~\ref{fig:S_one_half}.  This region is equal to the region $R$ in Figure~\ref{fig:BigPicture}. Let $\gamma$ be the dashed arc connecting the region $R$ to the $z$-pointed region.  Denote by $(a_{i,*})$, the intersection point between $a_{i,*}$ and $\gamma$.

\begin{lemma}\label{lem:admissibility}
The doubly-pointed Heegaard diagram described above, and appearing in Figures \ref{fig:S_one_half} and \ref{fig:BigPicture}, is weakly admissible.
\end{lemma}

\begin{proof}
Let $P$ be a nontrivial periodic domain for the pointed Heegaard diagram $(\Sigma, \beta, \alpha, w)$.  Suppose $P$ has nonzero local multiplicity in the $z$-pointed region.  Without loss of generality, we assume that this multiplicity is $+1$.  In particular, the multiplicity of the region just above the point $(a_{2,F})$ in Figure~\ref{fig:BigPicture} is $+1$.  In order for the $w$-pointed region to have multiplicity zero, the $\alpha_0$- and $\beta_0$-curves must be contained in the boundary of the periodic domain $P$ and must appear with multiplicity $\pm 1$ (depending on the chosen orientations of $\alpha_0$ and $\beta_0$).  This forces the small region southeast of $x_0$ in Figure~\ref{fig:S_one_half} to have multiplicity $+2$.  Since this region is the same as the region $R$ in Figure~\ref{fig:BigPicture}, $R$ must also have multiplicity $+2$.

Consider the dashed arc $\gamma$ connecting $R$ to the $z$-pointed region.  In order for $P$ to exist, the multiplicities of the regions intersected by $\gamma$ must go from $+2$ in the region $R$ to $+1$ in the $z$-pointed region.

However, the curve $\gamma$ intersects each $\alpha$-curve (other than $\alpha_0$) in two points, each with opposite sign.  The boundary of a periodic domain must be a sum of full $\alpha$- and $\beta$-curves, so if the local multiplicity of $P$ increases (or decreases) by a factor of $n$ as $\gamma$ passes one of the intersection points, it must decrease (or increase) by that same factor as $\gamma$ passes the other intersection point.  Therefore, the net change in the local multiplicity of the periodic domain $P$ along $\gamma$ between the region $R$ and the $z$-pointed region is zero.  We have seen that the multiplicity of the region $R$ is $+2$, whereas the multiplicity of the $z$-pointed region was assumed to be $+1$.  From this contradiction, we conclude that $P$ must have local multiplicity zero in the $z$-pointed region.

Since each $\alpha$- and each $\beta$-curve bound the $z$-pointed region on either side, and since $P$ has local multiplicity zero in the $z$-pointed region, we conclude that $P$ must have both positive and negative coefficients.
\end{proof}


\subsection{Computing $\STH(T)$} 
\label{sub:computing_T(T)}

Let $\phi \in \pi_2({\bf y}, {\bf y'})$ be a homotopy class with $n_z(\phi) = n_w(\phi) = 0$.  As is common in Heegaard Floer homology, we consider the ``shadow'' of $\phi$ on our Heegaard surface $\Sigma$.  This shadow is a sum of regions in the complement of the $\alpha$- and $\beta$-curves, and is denoted $D = \sum_i c_i D_i$.  If the homotopy class $\phi$ is to have a holomorphic representative $\mathbb{D}$, then each of the $c_i$ must be nonnegative.  In other words, $D$ is a positive domain and $\phi$ is a positive class.

Let ${\bf x} = (x_0, \dots, x_{2g})$; we show that the transverse invariant $\STH(T) = [(x_0, \dots, x_{2g})]$ is nonzero by proving that the generator ${\bf x}$ cannot appear in image of the Heegaard Floer differential.  This is accomplished by showing that the set of positive classes $\phi \in \pi_2({\bf y}, {\bf x})$ with $n_z(\phi) = n_w(\phi) = 0$ is empty for all ${\bf y} \neq {\bf x}$.

To draw a contradiction, we assume in what follows that $D$ is the domain of a positive class $\phi \in \pi_2({\bf y},{\bf x})$ for some generator ${\bf y} \neq {\bf x}$ with $n_z(\phi) = n_w(\phi) = 0$.

Suppose $R_1$ and $R_1$ are two adjacent regions in a Heegaard diagram (i.e. $R_1$ and $R_2$ share an edge), and $D$ is as above.  In general, the multiplicities of $R_1$ and $R_2$ can differ arbitrarily.  In our case, however, the multiplicities of any two adjacent regions can differ by at most one.  This is true because each of the $\alpha$- and $\beta$-curves in a Heegaard diagram coming from an open book decomposition bound the $z$-pointed region to either side.  Therefore, the boundary of any such region $D$ can never contain a full $\alpha$ or $\beta$ curve, and the multiplicities of $R_1$ and $R_2$ can differ by at most one.

Consider the region $R$ in Figure~\ref{fig:BigPicture}, and the curve $\gamma$ connecting $R$ to the $z$-pointed region.

\begin{lemma}\label{lem:nonnegative}
The net change in the local multiplicity of $D$ between the region $R$ and the $z$-pointed region along $\gamma$ is nonnegative.
\end{lemma}

\begin{proof}
The proof of Lemma \ref{lem:nonnegative} is similar to the proof of the admissibility lemma in Section \ref{sub:admissibility}.  	Recall that the curve $\gamma$ intersects each $\alpha_i$-curve in two points: $(a_{i,I})$ and $(a_{i,F})$.  We show that if the local multiplicity of the regions intersected by $\gamma$ decreases by a factor of $-1$ at the point $(a_{i,I})$, then there must be a corresponding increase in local multiplicity at the point $(a_{i,F})$.  Similarly, we show that if the local multiplicity of the regions intersected by $\gamma$ decrease by a factor of $-1$ at the point $(a_{i,F})$, then there must be a corresponding increase in local multiplicity at the point $(a_{i,I})$.

Observe that the local multiplicity along $\gamma$ cannot decrease as $\gamma$ passes over the point $(a_{1,I})$.  Since $D$ is the shadow of a positive class, if there is a decrease in local multiplicity at $(a_{1,I})$, a segment of the $\alpha_1$-curve between the intersection point $(a_{1,I})$ and $x_1$ must be contained in the boundary of the $D$.  Looking at the diagram in Figure~\ref{fig:BigPicture}, we see that any such arc has $z$-pointed region to the west, contradicting its existence.  

In the genus one case, a similar argument shows that there can be no decrease in the local multiplicity of $D$ at the point $(a_{2,I})$.  So assume that either the genus of $S$ is greater than one, or that we are considering an intersection point $(a)$ beyond $(a_{2,I})$ along $\gamma$.

Suppose that $(a) = (a_{i,I})$, and that the local multiplicity along $\gamma$ decreases by a factor $-1$ as it passes over at the point $(a)$.  Then, up to orientation, the segment of the $\alpha_i$-curve beginning at the point $(a)$ and traveling away from the center of the diagram to the point $x_i$ is contained in the boundary of $D$.  This implies that the region just past the intersection point $(a_{i,F})$ along $\gamma$ gains a $+1$ boost in local multiplicity.  Therefore, the increase in the local multiplicity at the point $(a_{i,F})$ balances the decrease in the local multiplicity at the point $(a_{i,I})$.

Similarly, if $(a) = (a_{i,F})$ and if the local multiplicity of $D$ along $\gamma$ decreases by a factor of $-1$ as $\gamma$ passes over $(a)$, then, again up to orientation, the segment of the $\alpha_i$-curve beginning at the point $(a_{i,I})$ and traveling away from the center of the diagram to the point $x_i$ is contained in the boundary of $D$.  This implies that the region just past the intersection point $(a_{i,I})$ gains a $+1$ boost in local multiplicity.  Thus, the decrease in local multiplicity at the point $(a_{i,F})$ is balanced by the increase in local multiplicity at the point $(a_{i,I})$.

Since each decrease in the local multiplicity of $D$ along $\gamma$ is balanced by a corresponding increase in local multiplicity somewhere else along $\gamma$, we have that the net change in the local multiplicity of $D$ between the region $R$ and the $z$-pointed region along the curve $\gamma$ is nonnegative.
\end{proof}

Consider the region $U$ in Figure~\ref{fig:BigPicture}, and the curve $\delta$ connecting $U$ to the $z$-pointed region.  By an argument similar to the proof of Lemma \ref{lem:nonnegative}, we have the following:

\begin{lemma}\label{lem:nonnegative2}
The net change in the local multiplicity of $D$ between the region $U$ and the $z$-pointed region along $\delta$ is nonnegative.
\end{lemma}
\begin{flushright}
$\Box$
\end{flushright}

\begin{proof}[Proof of \thmref{thm:binding}]
There are two main cases to consider.  

{\bf Case 1:}  Assume that the positive domain $D$ has nonzero local multiplicity in a region bordering the intersection point $x_0$.  In this case, the region $R$ in Figure~\ref{fig:BigPicture} has multiplicity $+1$.  By Lemma \ref{lem:nonnegative}, this implies that the multiplicity of the $z$-pointed region must be at least $+1$, a contradiction.

{\bf Case 2:}  Now suppose that $D$ has local multiplicity zero in all four of the regions bordering the point $x_0$.  In particular, this means that the multiplicity of the region $R$ is zero.  We investigate the possible configurations of $D$ near the center of Figure~\ref{fig:BigPicture}.

Suppose, for the moment, that all the regions bordered by the curve $\beta_0$ have zero multiplicity.  Then, near the center of Figure~\ref{fig:BigPicture}, the regions of $D$ with positive multiplicity are (locally) constrained to lie within the strip bounded by the darkened portions of the $\beta$-curves. 

In order for this to be the case, the boundary of the domain must have veered off the $\alpha$-curves while still contained within this strip.  Therefore, all the $\alpha$-curves are ``used up'' close to the center of the diagram (i.e. by the time they first intersect a darkened $\beta$-curve).  This, in turn, forces the multiplicity of the region $U$ in Figure~\ref{fig:BigPicture} to be positive.

By Lemma \ref{lem:nonnegative2}, this implies that the multiplicity of the $z$-pointed region must be positive, a contradiction.  Therefore, in order for such a nontrivial positive domain to exist, at least one of the regions bordered by $\beta_0$ must have nonzero multiplicity.

Recall that in Case 2 we are assuming that our domain $D$ is constant near $x_0$.  This means that the curve $\beta_0$ cannot be appear in the boundary of $D$ with nonzero multiplicity, so at least one of the regions intersected by $\gamma$ must have positive multiplicity.  Let $R'$ be the first region along $\gamma$ with positive multiplicity, and let $(a)$ be the $(a_{i,*})$ immediately preceding $R'$.

If $(a) = (a_{i,F})$, then by an argument similar to the proof of Lemma \ref{lem:nonnegative}, it can be shown that the net change in multiplicity between the region $R'$ and the $z$-pointed region must be nonnegative.  The fact that $(a)$ is a final point ensures that there can be no decrease in multiplicity at the point $(a_{i,I})$ since, by the definition of $R'$, the regions to both sides of this point have multiplicity zero.

On the other hand, suppose $(a) = (a_{i,I})$.  An argument similar to that in Lemma \ref{lem:nonnegative} demonstrates that for each decrease in multiplicity, there is a corresponding increase in multiplicity, except possibly at the point $(a_{i,F})$.  If the multiplicity decreases at the point $(a_{i,F})$, then the segment of the $\alpha_i$-curve from $(a_{i,F})$ to $(a)$ must be contained in the boundary of the domain.  This then implies that the multiplicity of the region $U$ is at least one.

Again, by Lemma \ref{lem:nonnegative2}, this forces the multiplicity of the $z$-pointed region to be positive, a contradiction.
\end{proof}

This completes the proof of \thmref{thm:binding}.


\section{The vanishing theorem} 
\label{sec:vanishing}

In this section, we prove \thmref{thm:vanishing}.  The proof in this case is similar to the proof of Theorem 4.5 in \cite{HKM2}. The key differences are that we must be careful to incorporate the Legendrian knot $L$ when choosing a Legendrian skeleton for the complement of the submanifold $N$, and that we must be cautious about the changes made to the diagram in the spinning process used to make the diagram strongly admissible.

\begin{proof}[Proof of \thmref{thm:vanishing}]
We begin by constructing a partial open book decomposition for the contact submanifold $(N,\xi|_N)$, which can be extended to an open book decomposition for all of $(Y,\xi)$.  Following \cite{HKM2}, we must show that the basis $\{a_1, \dots, a_r\}$ for the partial open book decomposition of $(N,\xi|_N)$ can be extended to a basis $\{a_1, \dots, a_r, a_0', a_1', \dots, a_s'\}$ for the extended open book decomposition of $(Y,\xi)$, where $L \cap (\bigcup a_i \cup \bigcup a_j') = L \cap a_0' = 1 \mr{ pt}$.\\

{\bf Claim:}  We may assume without loss of generality that the complement of $N$ is connected.

\begin{proof}[Proof of claim]
Let $(M,\xi)$ be a compact manifold with possibly nonempty boundary, and let $(M',\xi|_{M'})$ be a compact submanifold of $(M,\xi)$ with convex boundary.  In \cite{HKM3}, the authors show that the vanishing of the contact invariant for $(M',\xi|_{M'})$ implies the vanishing of the contact invariant for $(M,\xi)$.  

Suppose the complement of $N$ is disconnected.  Then, since $c(N,\xi|_N) = 0$, the contact manifold obtained by gluing the components of $Y - N$ not containing $L$ to $N$ must also have vanishing contact invariant.  In particular, we may assume without loss of generality that $Y - N$ is connected.
\end{proof}

\begin{figure}[htbp]
	\centering
	\psfrag{a}{$L$}
	\psfrag{b}{$\nu(K)$}
	\psfrag{c}{\hspace{-10px}$\nu(K')$}
	\psfrag{d}{$N$}
	\psfrag{f}{$\partial N$}
		\includegraphics[scale=0.8]{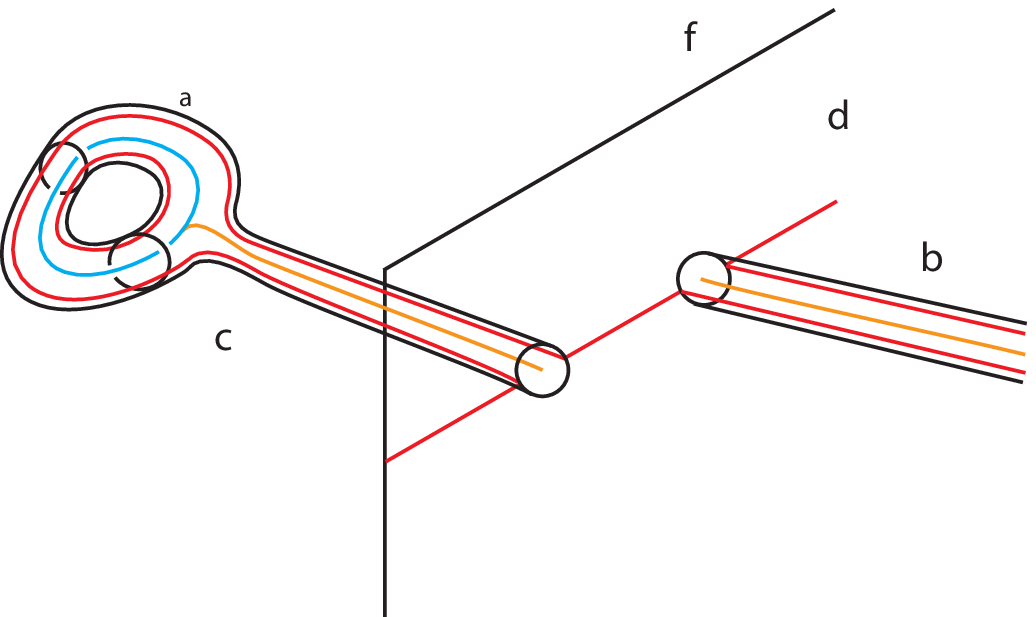}
		\caption{}
		\label{fig:skeleton}
\end{figure}

Let $K$ be a Legendrian skeleton for $N$, and let $K'$ be an extension of the Legendrian knot $L$ to a Legendrian skeleton for $N' = Y - N$ (see Figure~\ref{fig:skeleton}).  Assume that the univalent vertices of $K$ and $K'$ in $\partial N$ do not intersect.

The Legendrian skeleton $K$ gives us a partial open book decomposition for $(N, \xi|_N)$.  Let $\nu(K)$ be a standard neighborhood of $K$ inside of $N$, and let $\nu(K')$ be a standard neighborhood of $K'$ inside of $Y - N$.  We can build an open book decomposition for all of $Y$ by considering the the handlebodies $(N' - \nu(K')) \cup \nu(K)$ and $\nu(K') \cup (N - \nu(K))$.  By construction, each of these handlebodies are disk decomposable.  A page $S$ of the open book for $(Y,\xi)$ is constructed from the page of the partial open book for $(N,\xi|_N)$ by repeatedly attaching 1-handles away from the portions of the open book coming from the boundary of $\nu(K)$.  This construction is depicted in Figure~\ref{fig:page}.

In Figure~\ref{fig:page}, the portion of the page of the open book coming from the boundary of $N$ is shown in black, and has its boundary lines thickened.  The portion of the page coming from the boundary of $\nu(K)$ is lightly colored (orange), and appears in the lower right portion of the figure.  Finally, the portion of the page coming from the extension of the open book to all of $Y$ is also lightly colored (green), and appears in the lower left corner of the figure.  

Let $\{a_1,\dots,a_n\}$ be a basis for the partial open book coming from $(N,\xi|_N)$, and let $\phi$ be the corresponding partially defined monodromy map for this open book.  Consider a new partial open book, whose page is equal to $S$, and whose partially defined monodromy map is equal to $\phi$.  Because this new partial open book only differs from the partial open book coming from $(N,\xi|_N)$ by handle attachments away from $\partial \nu(K)$, the contact element for this new partial open book vanishes along with $c(N,\xi|_N)$.

Since $Y - N$ is connected, the basis $\{a_1, \dots, a_n \}$ can, after a suitable number of stabilizations, be extended to a basis for all of $Y$.  

\begin{figure}[htbp]
	\centering
	\psfrag{a}{$\partial N$}
	\psfrag{b}{\hspace{-15px}$\partial \nu(K')$}
	\psfrag{c}{$\partial \nu(K)$}
	\psfrag{e}{\hspace{-2px}$a_0'$}
	\psfrag{g}{$L$}
	\psfrag{d}{$S$}
	\includegraphics[scale=0.85]{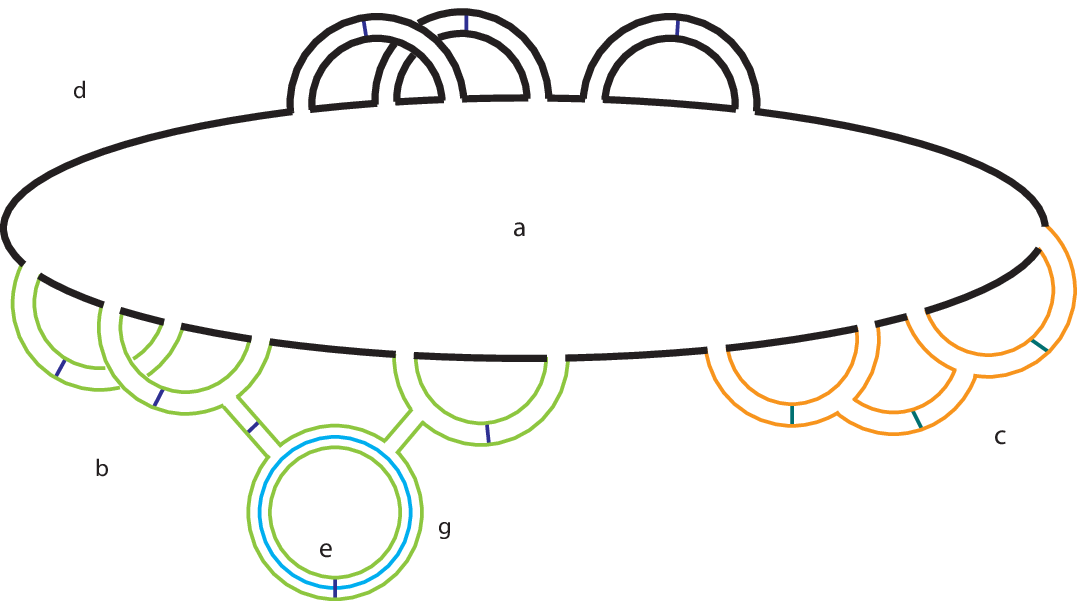}
	\caption{}
	\label{fig:page}
\end{figure}

By construction, the new monodromy map $\phi'$ extends $\phi$, the monodromy map for $N$.  We can see our Legendrian $L$ on the page $S$.  The local picture around $L \subset S$ (shown in blue) must look like that in Figure~\ref{fig:page}.

As was observed in \cite{LOSS}, the ``spinning'' isotopies needed to make this Heegaard diagram strongly admissible can be performed on the portion of the Heegaard diagram coming from the page $S_1$.  This changes the monodromy map $\phi'$, but only within its isotopy class.

If we delete the $\alpha$- and $\beta$-curves coming from $\{a_0',a_1', \dots, a_s'\}$, then we are left with a diagram which is essentially equivalent to that coming from the partial open book $(S,\phi)$, but whose monodromy has been changed by an isotopy.  Since altering the monodromy map by an isotopy cannot change whether or not the contact element vanishes in sutured Floer homology, we know that the contact element corresponding to the partial open book $(S,\phi')$ vanishes.  That is, if ${\bf x} = (x_1, \dots, x_n)$, then there exist $c_i$ and ${\bf y}_i$ such that $\partial (\sum_i c_i {\bf y}_i) = {\bf x}$ in the sutured Floer homology of the manifold obtained from the partial open book $(S, \phi')$.  

Let ${\bf x}' = (x_0', x_1', \dots, x_s')$; we claim that $\partial (\sum_i c_i ({\bf y}_i, {\bf x}')) = ({\bf x}, {\bf x}')$ in $\HFKM(-Y,L)$.  The intersection points coming from ${\bf x}'$ must map to themselves via the constant map.  This allows us to ignore the $\alpha$- and $\beta$-curves corresponding to these intersection points, leaving us with a diagram which is essentially equivalent to the partial open book $(S,\phi')$.
\end{proof}

\bibliographystyle{amsalpha}
\nocite{*}
\bibliography{TransInvt}

\providecommand{\bysame}{\leavevmode\hbox to3em{\hrulefill}\thinspace}
\providecommand{\MR}{\relax\ifhmode\unskip\space\fi MR }
\providecommand{\MRhref}[2]{%
  \href{http://www.ams.org/mathscinet-getitem?mr=#1}{#2}
}
\providecommand{\href}[2]{#2}
\begin{thebibliography}{HKM09b}

\bibitem[Col08]{Co}
Vincent Colin, \emph{Livres ouverts en g\'eom\'etrie de contact (d'apr\`es
  {E}mmanuel {G}iroux)}, Ast\'erisque (2008), no.~317, Exp. No. 969, vii,
  91--117, S{\'e}minaire Bourbaki. Vol. 2006/2007.

\bibitem[EFM01]{EFM}
Judith Epstein, Dmitry Fuchs, and Maike Meyer, \emph{Chekanov-{E}liashberg
  invariants and transverse approximations of {L}egendrian knots}, Pacific J.
  Math. \textbf{201} (2001), no.~1, 89--106.

\bibitem[EH01]{EH}
John~B. Etnyre and Ko~Honda, \emph{Knots and contact geometry. {I}. {T}orus
  knots and the figure eight knot}, J. Symplectic Geom. \textbf{1} (2001),
  no.~1, 63--120.

\bibitem[ENV10]{ENV}
John~B. Etnyre, Lenhard~L. Ng, and Vera V{\'e}rtesi, \emph{Legendrian and
  transverse twist knots}, Preprint, \href{http://arxiv.org/abs/1002.2400}{\tt
  arXiv:1002.2400} [math.SG], 2010.

\bibitem[Etn05]{Et2}
John~B. Etnyre, \emph{Legendrian and transversal knots}, Handbook of knot
  theory, Elsevier B. V., Amsterdam, 2005, pp.~105--185.

\bibitem[Etn06]{Et1}
\bysame, \emph{Lectures on open book decompositions and contact structures},
  Floer homology, gauge theory, and low-dimensional topology, Clay Math. Proc.,
  vol.~5, Amer. Math. Soc., Providence, RI, 2006, pp.~103--141.

\bibitem[EV10]{EV}
John~B. Etnyre and David~Shea {Vela--Vick}, \emph{Torsion and open book
  decompositions}, Internat. Math. Res. Notices (2010).

\bibitem[EVZ10]{EVZ}
John~B. Etnyre, David~Shea {Vela--Vick}, and Rumen Zarev, \emph{Sutured
  {L}egendrian invariants and invariants of open contact manifolds}, In
  Preparation, 2010.

\bibitem[GHV07]{GHV}
Paolo Ghiggini, Ko~Honda, and Jeremy {Van Horn--Morris}, \emph{The vanishing of
  the contact invariant in the presence of torsion}, Preprint,
  \href{http://arxiv.org/abs/0706.1602v2}{\tt arXiv:0706.1602v2} [math.GT],
  2007.

\bibitem[Gir02]{Gi}
Emmanuel Giroux, \emph{G\'eom\'etrie de contact: de la dimension trois vers les
  dimensions sup\'erieures}, Proceedings of the International Congress of
  Mathematicians, Vol. II (Beijing, 2002) (Beijing), Higher Ed. Press, 2002,
  pp.~405--414.

\bibitem[HKM08]{HKM3}
Ko~Honda, William Kazez, and Gordana Mati\'c, \emph{Contact structures, sutured
  {F}loer homology, and {TQFT}}, Preprint,
  \href{http://arxiv.org/abs/0807.2431}{\tt arXiv:0807.2431} [math.GT], 2008.

\bibitem[HKM09a]{HKM1}
\bysame, \emph{On the contact class in {H}eegaard {F}loer homology}, J.
  Differential Geom. \textbf{83} (2009), no.~2, 289--311.

\bibitem[HKM09b]{HKM2}
\bysame, \emph{The contact invariant
  in sutured {F}loer homology}, Invent. Math. \textbf{176} (2009), no.~3,
  637--676.

\bibitem[Lip06]{Li}
Robert Lipshitz, \emph{A cylindrical reformulation of {H}eegaard {F}loer
  homology}, Geom. Topol. \textbf{10} (2006), 955--1097 (electronic).

\bibitem[LOSS09]{LOSS}
Paolo Lisca, Peter Ozsv{\'a}th, Andr{\'a}s~I. Stipsicz, and Zolt{\'a}n
  Szab{\'o}, \emph{Heegaard {F}loer invariants of {L}egendrian knots in contact
  three-manifolds}, J. Eur. Math. Soc. (JEMS) \textbf{11} (2009), no.~6,
  1307--1363.

\bibitem[OS04a]{OS3}
Peter Ozsv{\'a}th and Zolt{\'a}n Szab{\'o}, \emph{Holomorphic disks and knot
  invariants}, Adv. Math. \textbf{186} (2004), no.~1, 58--116.

\bibitem[OS04b]{OS2}
\bysame, \emph{Holomorphic disks and three-manifold invariants: properties and
  applications}, Ann. of Math. (2) \textbf{159} (2004), no.~3, 1159--1245.

\bibitem[OS04c]{OS1}
\bysame, \emph{Holomorphic disks and topological invariants for closed
  three-manifolds}, Ann. of Math. (2) \textbf{159} (2004), no.~3, 1027--1158.

\bibitem[OS05]{OS4}
\bysame, \emph{Heegaard {F}loer homology and contact structures}, Duke Math. J.
  \textbf{129} (2005), no.~1, 39--61.

\bibitem[OS08]{OSt}
Peter Ozsv{\'a}th and Andr{\'a}s~I. Stipsicz, \emph{Contact surgeries and the
  transverse invariant in knot {F}loer homology}, Preprint,
  \href{http://arxiv.org/abs/0803.1252}{\tt arXiv:0803.1252v1} [math.SG], 2008.

\bibitem[Ras03]{Ra}
Jacob Rasmussen, \emph{Floer homology and knot complements}, Ph.D. thesis,
  Harvard University, 2003, \href{http://arxiv.org/abs/math/0306378}{\tt
  arXiv:math/0306378v1} [math.GT].

\bibitem[SV09]{StV}
Andr{\'a}s~I. Stipsicz and Vera V{\'e}rtesi, \emph{On invariants for
  {L}egendrian knots}, Pacific J. Math. \textbf{239} (2009), no.~1, 157--177.

\bibitem[Tor00]{To}
Ichiro Torisu, \emph{Convex contact structures and fibered links in
  3-manifolds}, Internat. Math. Res. Notices (2000), no.~9, 441--454.

\bibitem[TW75]{TW}
William~P. Thurston and Horst~E. Winkelnkemper, \emph{On the existence of
  contact forms}, Proc. Amer. Math. Soc. \textbf{52} (1975), 345--347.

\end{thebibliography}


\end{document}